\documentclass[11pt, oneside]{article}   	
\usepackage{geometry}                		
\usepackage{graphicx}				
\usepackage{amssymb}
\usepackage{amsmath}
\usepackage{subfig}
\usepackage{caption}
\usepackage{comment}

\usepackage[utf8]{inputenc}
\usepackage{amsmath}
\usepackage{graphicx}
\usepackage{float}	
\usepackage{amsmath}

\numberwithin{equation}{section}
\begin{document}
\title{Two Matrix Model, the Riemann Hypothesis \\ and Master Matrix Obstruction\\
}
\author{ Michael McGuigan\\
email contact: michael.d.mcguigan@gmail.com
}

\date{}
\maketitle
\begin{abstract}
We identify the Riemann Xi function as the Baker-Akhiezer function for a $(p,1)$ two matrix model as $p$ goes to infinity. We solve the two matrix model using biorthogonal polynomials and study the zeros of the polynomials in the double scaling limit as $N$ goes to infinity. We find zeros off the critical line at finite $N$ which possibly go to infinity as $N$ goes to infinity. We study other Baker-Akhiezer functions whose zeros are known to be on a critical line using the two matrix model technique and find the zeros on the critical line in those cases. We study other L-functions using the two matrix model and compare the biorthogonal method with other approaches to the two matrix model such as the master matrix approach and saddle point method. In cases where there are zeros off the critical line  the master matrix approach encounters an obstruction to the solution to a quenched master matrix.
\end{abstract}
\newpage

\section{Introduction}

The Riemann hypothesis, which states that the  zeros of the Riemann zeta function lie on the critical line, continues to intrigue and attract mathematicians and physicists alike \cite{Conrey}
\cite{Konig}
\cite{Medlin}
\cite{Cardon}
\cite{Perez-Moure}
\cite{Conrey2}. There are indications that the hypothesis is related to  quantum mechanics, although it appeared in the  the nineteenth century, well before quantum mechanics had been formulated. In addition the spacing of the zeros have distributions in common with random matrix integrals. In this paper we explore a different connection of the Riemann hypothesis and quantum mechanics  through its relation to the $(p,1)$ two matrix model in the limit where $p$ and the size of the matrix $N$ goes to infinity. Here $p$ refers to the degree of the potential of one of the matrices in the two matrix model with the other matrix having a degree $1$ or linear potential. The relation will allow us to relate the Riemann Xi function to the Baker-Akhiezer of the two matrix model and study its zeros using the biorthogonal polynomial method of solution.

This paper is organized as follows. In section 2 we review some basic facts about the Riemann Xi function and its representation as a Fourier integral. In section 3 we review some features of the two matrix model needed for this paper and how the representation of the Baker-Akhiezer function as a Fourier integral leads to a connection to the Riemann Xi function \cite{McGuigan:2007pr}
\cite{McGuiganMaster}. In section 4 we apply the method of biorthogonal polynomials \cite{biorthogonal1}
\cite{biorthogonal2} to the study of the simplest case, the $(2,1)$ two matrix model, whose Baker-Akhiezer function is the Airy function and show how to obtain the zeros as the zeros of a polynomial. In section 5 we apply the same method to the Riemann Xi function, except in that case we find zeros off the critical line which possibly go to infinity as $p$ and $N$ go to infinity. In section 6 we apply the method to other L functions such as the Ramanujan L function and also see zeros off the critical line  which also may go to infinity in the double scaling limit. In section 7 we discuss Baker-Akhiezer function which are know to have zeros on critical lines and show in those cases the two matrix model also yields zeros on the critical line. In section 8 we discuss the $p$ goes to infinity limit in a case where the polynomial expansion of the potential is particularly straightforward again finding zeros on the corresponding critical line. In section 9 we discuss other methods to solve the two matrix model besides the biorthogonal polynomial method such as the master matrix and saddle point method. In the case of the master matrix method we find an obstruction to the solution of the master matrix when some zeros are off the critical line. Finally in section 10 we discuss our conclusions and directions for future work.

\section{Review of Riemann Xi function}

The Riemann Xi function is defined as 
\begin{equation}
\Xi (z) = \xi (\frac{1}{2} + iz)
\end{equation}
where
\begin{equation}
\xi (s) = \frac{1}{2}s(s - 1){\pi ^{ - s/2}}\Gamma (s/2)\zeta (s)
\end{equation}
It shares all the nontrivial zeros of the Riemann zeta function but does not have the pole or the trivial zeros. It can be represented as a Fourier integral as:
\begin{equation}
\psi_R(z)=\Xi (z) = \int_{ - \infty }^\infty  {\Phi (x){e^{ixz}}dx} 
\end{equation}
where
\begin{equation}
\Phi (x) = \sum\limits_{n = 1}^\infty  {\left( {4{\pi ^2}{n^4}{e^{9x/2}} - 6\pi {n^2}{e^{5x/2}}} \right)} \exp \left( { - \pi {n^2}{e^{2x}}} \right)
\end{equation}
with derivative:
\begin{equation}
\Phi '(x) = \sum\limits_{n = 1}^\infty  {\left( {30{\pi ^2}{n^4}{e^{9x/2}} - 15\pi {n^2}{e^{5x/2}} - 8{\pi ^3}{n^6}{e^{13x/2}}} \right)} \exp \left( { - \pi {n^2}{e^{2x}}} \right)
\end{equation}
For connections to the $(p,1)$ two matrix model one can also form the series:
\begin{equation}
U(x) =  - \log \left( {\Phi (x)} \right) = \sum\limits_{n = 0}^\infty  {{a_n}{x^n}} 
\end{equation}
and the truncated version at order $p+1$ given by:
\begin{equation}
{
U_p}(x) = \sum\limits_{n = 0}^{p + 1} {{a_n}{x^n}} \end{equation}
It is convenient to rescale the term of order $p+1$ so its coefficient is given by $1/(p+1)$. Then we have:
\begin{equation}
{U_p}\left( {\frac{{\bar x}}{{{{\left( {{a_{p + 1}}\left( {p + 1} \right)} \right)}^{1/(p + 1)}}}}} \right) = \frac{{{{\bar x}^{p + 1}}}}{{p + 1}} + \sum\limits_{n = 2}^{p-1} {{s_{n - 1}}{{\frac{{\bar x}}{n}}^n}} +a_0 
\end{equation}
where:
\begin{equation}
{s_{n - 1}} = n{a_n}\frac{1}{{{{\left( {{a_{p + 1}}\left( {p + 1} \right)} \right)}^{n/(p + 1)}}}}
\end{equation}
The Baker-Akhiezer representation of the Riemann Xi function is then:
\begin{equation}
{\psi _R}(z) = \Xi(z) = {\lim _{p \to \infty }}\int_{ - \infty }^\infty  {{e^{ - {U_p}(x)}}{e^{izx}}dx} 
\end{equation}
which as we shall see is the limit of the expectation value of the characteristic polynomial for the $(p,1)$ two matrix model as $p$ goes to infinity.

\section{Review $(p,1)$ two matrix model}

We will mainly follow the approach to the $(p,1$ two matrix model in \cite{Hashimoto:2005bf}
\cite{Maldacena:2004sn}. Many aspects two matrix model have been developed in the physics [13-77] and mathematics literature [78-84].
The partition function for the $(p,1)$ two matrix model is given by
\begin{equation}Z = \int {dAdB{e^{ - \frac{1}{g}Tr(V(A )+ B - AB)}}} \end{equation}
where $A$ and$B$ are two Hermitian $N \times N$ matrices,
$V(A)$ is the potential for the $A$ matrix and $B$ is the potential for the $B$ matrix. After shifting $A$ by the identity the partition function becomes:
\begin{equation}Z = \int {dAdB{e^{ - \frac{1}{g}Tr(V(A + I) - AB)}}} \end{equation}
The quantity that interests us is the expectation value of the characteristic polynomial of the $B$ matrix given by:

\begin{equation}
\left\langle {\det (b - B)} \right\rangle  = \int {dAdB{e^{ - \frac{1}{g}Tr(V(A + I) - AB)}}\det (b - B)} 
\end{equation}
The expectation value of the characteristic polynomial of $B$ can be computed using biorthogonal polynomials $P_m(a),Q_n(b)$ which obey
\begin{equation}
\int {dadb{P_m}} (a){Q_n}(b){e^{ - \frac{1}{g}\left( {V(a + 1) - ab)} \right)}} = {h_m}{\delta _{m,n}}
\end{equation}
The biorthogonal  polynomials are given as:
\begin{align}
&{P_m}(a) = {a^m}\nonumber \\
&{Q_n}(b) = {\left( { - g\frac{\partial }{{\partial a}}} \right)^n}{\left. {\left( {{e^{ - \frac{1}{g}\left( {V(a + 1) - ab)} \right)}}} \right)} \right|_{a = 0}}
\end{align}
And the expectation value of the characteristic polynomial can be shown to be:
\begin{equation}
\left\langle {\det (b - B)} \right\rangle  = {Q_N}(b) = {\left( { - g\frac{\partial }{{\partial a}}} \right)^N}{\left. {\left( {{e^{ - \frac{1}{g}\left( {V(a + 1) - ab)} \right)}}} \right)} \right|_{a = 0}}
\end{equation}
To obtain the Baker-Akhiezer function we take the double scaling limit where the potential $V(A)$ is a $p$ degree matrix polynomial of the form:
\begin{equation}
V(A) = {V_p}(A) + {s_1}{\epsilon ^{p - 1}}{V_1}(A) + {s_3}{\epsilon ^{p - 2}}{V_3}(A) +  \ldots  + {s_{p - 2}}{\epsilon ^2}{V_{p - 2}}(A)
\end{equation}
where:
\begin{equation}
{V_n}(A) =  - \sum\limits_{k = 1}^n {\left( {\frac{{{A^k}}}{k} - \frac{I}{k}} \right)} 
\end{equation}
and 
\begin{align}
&\epsilon  = \left( \frac{1}{N}\right)^{ \frac{1}{p + 1}}\nonumber \\
&g = \frac{1}{N} + \frac{1}{N}\left( {{\epsilon ^2}{s_{p - 2}} + {\epsilon ^4}{s_{p - 4}} +  \ldots  + {\epsilon ^{p - 1}}{s_1}} \right)\nonumber \\
&a =  - 1 + \epsilon x
\end{align}
Then taking the large $N$ limit as $\epsilon$ goes to zero one can express the expectation value of the characteristic polynomial as the Fourier integral
\begin{equation}
\psi (z) = \int_{ - \infty }^\infty  {{e^{ - {U_p}(x)}}{e^{izx}}dx} 
\end{equation}
with $U_p(x)$ given by:
\begin{equation}
{U_p}(x) = \frac{x^{p+1}}{p+1} + {s_1}\frac{{{x^2}}}{2} + {s_1}\frac{{{x^4}}}{4} \ldots  + {s_{p - 2}}\frac{{{x^{p - 1}}}}{{p - 1}}
\end{equation}
which represents the Baker-Akhierzer function as a Fourier integral.  In the following sections we will apply these formulas to the study of various  functions such as the Riemann Xi function which be represented as the Baker-Akhiezer function of two matrix models.

If one defines the Jacobi Matrix $J_M$
\begin{equation}
x\frac{{{Q_n}(x)}}{{\sqrt {{h_n}} }} = \sum\limits_{m = 0}^{n + 1} {{(J_M)_{n,m}}\frac{{{Q_m}(x)}}{{\sqrt {{h_m}} }}} \end{equation}
Then taking the $N\times N$ submatrix of $J_M$ yields:
\begin{equation}
\left\langle {\det (b - B)} \right\rangle  = \det(b-J_M)
\end{equation}
one can also obtain the polynomials $Q_n(b)$ from the generating function:
\begin{equation}
{e^{\frac{1}{g}V(gt + 1)}}{e^{ - yt}} = \sum\limits_{N = 1}^\infty  {\frac{{{Q_N}(y){t^N}}}{{N!}}} 
\end{equation}

\section{Two Matrix Model and Airy function zeros}

The $(2,1)$ Matrix model was studied in \cite{Maldacena:2004sn}. Because $p=2$  the derivative formula for the $Q_n(b)$ polynomials reduces to the Hermite polynomials and the expectation value of the characteristic polynomial is given by:
\begin{equation}
\left\langle {\det \left( {b - B} \right)} \right\rangle  = {\left( {\frac{g}{4}} \right)^{N/2}}{H_N}\left( {b\frac{1}{{\sqrt g }}} \right)
\end{equation}
For $N=16$ and $gN=1$ this leads to the polynomial:
\begin{align}
Q_{16}(b) = b^{16}-3.75 b^{14}+5.33203 b^{12}-&3.66577 b^{10} +1.28875 b^8-0.225531 b^6 \nonumber\\
&+0.0176196 b^4-0.000471954 b^2+1.84357 \times 10^{-6} 
\end{align}
The zeros of the polynomial are:
\begin{align}
&b_1 = -1.17218 \nonumber \\
&b_2 = -0.967362 \nonumber \\ 
&b_3 = -0.79425 \nonumber \\
&b_4 = -0.636551 \nonumber \\
&b_5 = -0.487947 \nonumber \\
&b_6 = -0.345065 \nonumber \\
&b_7 = -0.205738 \nonumber \\
&b_8 = -0.0683703 \nonumber \\
&b_9 = 0.0683703 \nonumber \\     
&b_{10} = 0.205738 \nonumber \\
&b_{11} = 0.345065 \nonumber \\
&b_{12} = 0.487947 \nonumber \\
&b_{13} = 0.636551 \nonumber \\
&b_{14} = 0.79425 \nonumber \\
&b_{15} = 0.967362 \nonumber \\
&b_{16} = 1.17218 \\
\end{align}
The Baker-Akhiezer function for the $(2,1)$ two matrix model is the Airy function. The zeros of the polynomial are related to the Airy function $Ai(y)$ through a linear relation $y=Ab+c$. For 
\begin{equation}
y = -8 (2^{1/6}) \left(\sqrt{2}-b\right)
\end{equation}
the first three  zeros are:
\begin{align}
&y_1 =-8  (2^{1/6})\left(\sqrt{2}-b_{16}\right) \nonumber \\
&y_2 =-8 (2^{1/6}) \left(\sqrt{2}-b_{15}\right) \nonumber \\
&y_3 =-8 (2^{1/6}) \left(\sqrt{2}-b_{14}\right)
\end{align}
yielding:
\begin{align}
&y_1 =-2.17335 \nonumber \\
&y_2=-4.01259 \nonumber \\
&y_3 = -5.56709
\end{align}
which can be compared to the exact value for Airy zeros given by -2.33811, -4.08795, -5.52056. One can obtain more accurate values of the Airy zeros from the polynomials by going to larger values of N as shown in figure 1 for $N=32$.
\begin{figure}[!htb]
\centering
\minipage{0.5\textwidth}
  \includegraphics[width=\linewidth]{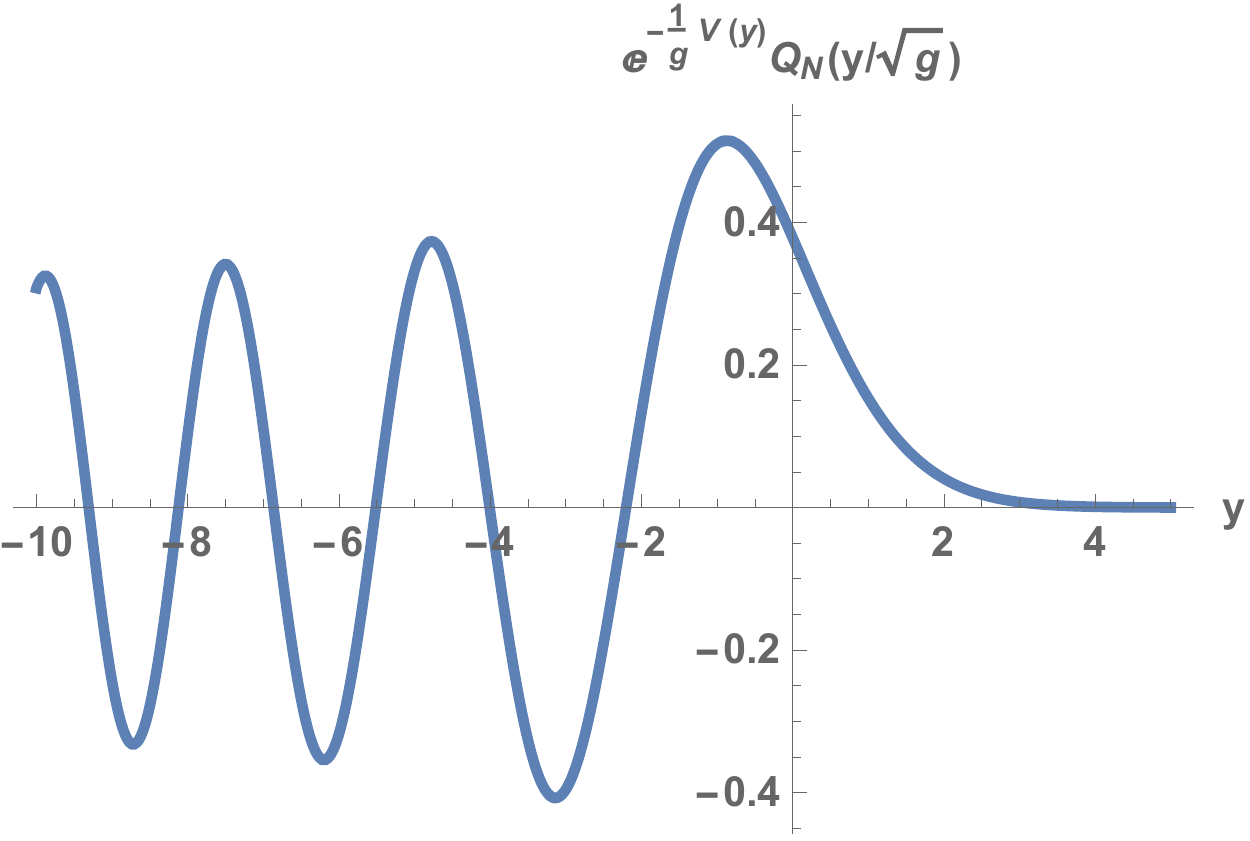}
\endminipage\hfill
\minipage{0.48\textwidth}
  \includegraphics[width=\linewidth]{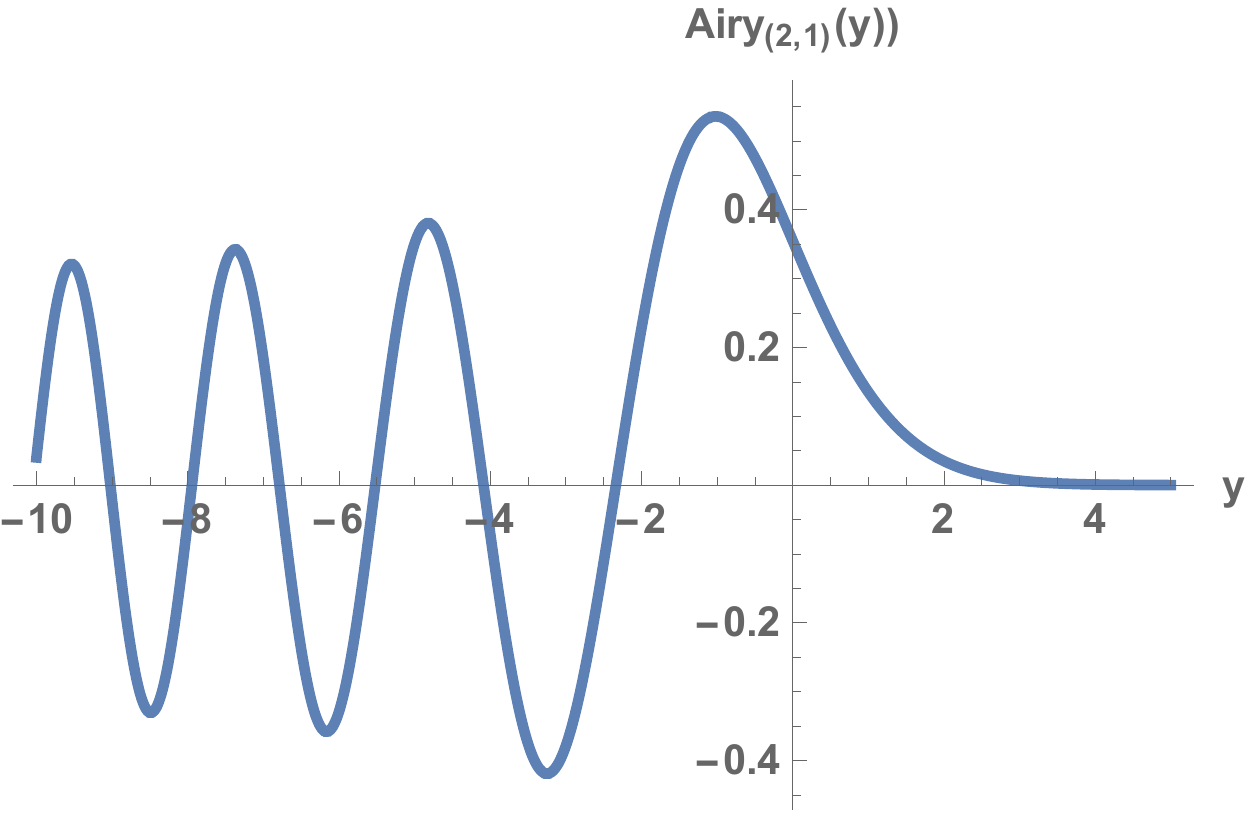}
\endminipage\hfill
\caption{(left) Approximation to the Airy function from the $(2,1)$ two Matrix model with $N=32$. (right) Exact expression for the Airy function. The two Matrix model approximation gives a good estimate to the first five zeros for $N=32$. 
}
\end{figure}

\section{ Two Matrix Model and Riemann Xi function zeros}

The logarithm of the Phi function has the expansion
\begin{equation}
- \log (\Phi (x)) = \sum\limits_{n = 0}^\infty  {{a_n}{x^n}}  = {a_0} + {a_2}{x^2} + {a_4}{x^4} + {a_6}{x^6} + {a_8}{x^8} +  \ldots 
\end{equation}
working with the sum at the eighth order we have:
\begin{align}
&a_0=0.112728 \nonumber \\
&a_2= 9.3634 \nonumber \\
&a_4=5.95896 \nonumber \\
&a_6= -2.09194 \nonumber \\
&a_8=3.53296 
\end{align}
then using the expansion:
\begin{equation}
{U_p}(\bar x) = \frac{{{{\bar x}^{p + 1}}}}{{p + 1}} + {s_1}\frac{{{{\bar x}^2}}}{2} + {s_3}\frac{{{{\bar x}^4}}}{2} \ldots  + {s_{p - 1}}\frac{{{{\bar x}^{p - 1}}}}{2}
\end{equation}
for $p=7$ we obtain the $s_n$ coefficients for the corresponding $(7,1)$ two matrix model as:
\begin{align}
&s_1=8.12192\nonumber \\
&s_4= 4.48349\nonumber \\
&s_6=-1.02395
\end{align}
Then using formulas (3.5) and (3.9) we have for the $N=16$ Q polynomial for $p=7$:

\begin{align}
Q_{16}[b]= b^{16}&+141.088 b^{15}+8952.1 b^{14}+338149. b^{13}+8.48406\times 10^6 b^{12} \nonumber \\
&+ 1.49383\times 10^8 b^{11}+1.90155\times 10^9 b^{10} +1.77654\times 10^{10} b^9 +1.2243\times 10^{11} b^8 \nonumber \\
&+ 6.20423\times 10^{11} b^7+2.28714\times 10^{12} b^6+6.01787\times 10^{12} b^5 +1.09783\times 10^{13} b^4 \nonumber\\
&+ 1.33068\times 10^{13} b^3+1.00497\times 10^{13} b^2+4.23563\times 10^{12} b+7.61563\times 10^{11}
\end{align}
The polynomial zeros are given by: 
\begin{align}   
&b_1 = -22.8352 \nonumber \\
&b_2 = -19.7169 \nonumber\\
&b_3 = -17.2138 \nonumber\\
&b_4 = -15.0434 \nonumber\\
&b_5 = -13.0978 \nonumber\\
&b_6 = -11.322 \nonumber\\
&b_7 = -9.68381 \nonumber\\
&b_8 = -8.16324 \nonumber\\
&b_9 = -6.74774 \nonumber\\
&b_{10} = -5.43011 \nonumber\\
&b_{11} = -4.20771 \nonumber\\
&b_{12}  = -3.08255 \nonumber\\
&b_{13} = -2.0612 \nonumber\\
&b_{14} = -1.12636 \nonumber\\
&b_{15} = -0.677917-0.213125 i \nonumber\\
&b_{16} = -0.677917+0.213125 i 
\end{align}
Fitting the first two zeros to the Riemann zeta zeros of the form  $z_i=Ab_i+c$ we find $A=2.20867$ and $c=64.5702$ then the first three zeros determines by the polynomial zeros are $14.1347, 21.022, 26.5505$ which can be compared with the exact values $14.1347, 21.022, 25.0109 $. One can increase the accuracy of the polynomial zero estimation of the zeros by going to larger $p$ and higher values of $N$. A plot of the Riemann Xi function which is the Baker Akheizer function  for the two matrix model is shown in figure 2.

%
\begin{figure}
\centering
  \includegraphics[width = .75 \linewidth]{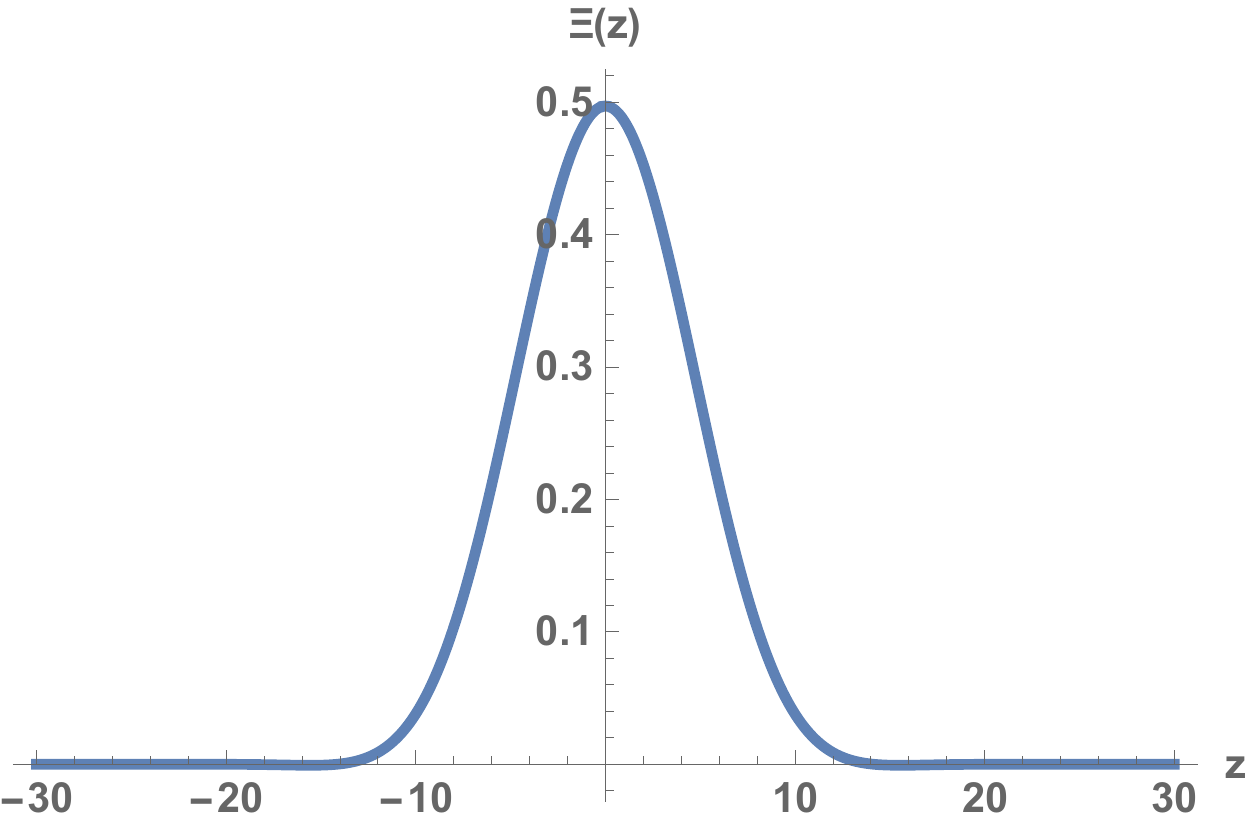}
  \caption{Plot of the Riemann Xi function which is the Baker Akheizer function  for a two matrix model associated with the Riemann zeta function. }
  \label{fig:Radion Potential}
\end{figure}
\section{ Two Matrix model and other L-function zeros}

We can use the same two Matrix model technique to study other L-functions. For example the Ramanujan tau L function can be represented as a Fourier integral through:
\begin{equation}\psi_L(z)= {\Xi _L}(z) = {(2\pi )^{ - (6 + iz)}}\Gamma (6 + iz)L(6 + iz) = \int_{ - \infty }^\infty  {{\Phi _L}(x){e^{ixz}}dx} \end{equation}
with $L$ the Ramanujan tau L function and:
\begin{equation}{\Phi _L}(x) = {e^{ - 6x}}\Delta (i{e^{ - x}}) = {e^{ - 6x}}{e^{ - 2\pi {e^{ - x}}}}\prod\limits_{n = 1}^\infty  {{{\left( {1 - {e^{ - 2\pi n{e^{ - x}}}}} \right)}^{24}}} \end{equation}
Defining:
\begin{equation}
    {\Phi _L}(x) = e^{-U_L(x)}
\end{equation}
we find $U_L$ has an expansion out to eighth order given by :
\begin{equation}
U_L^{(8)}(x) = 0.112629 x^8-0.188291 x^6+0.971962 x^4+3.89463 x^2+6.32813
\end{equation}
so:
\begin{align}
    &a_0 = 6.32813 \nonumber \\
    &a_2 = 3.89463 \nonumber \\
    &a_4 = 0.971962 \nonumber \\
    &a_6 = -0.188291 \nonumber \\
    &a_8 = 0.112629
\end{align}
and these determine  the parameters of the $(7,1)$ two Matrix model to be:
\begin{align}
    &s_1 = 7.99487 \nonumber \\
    &s_3 = 4.0958 \nonumber \\
    &s_5 = -1.22159 
\end{align}
Using formula (3.5) for $N=16$ we obtain the polynomial 
\begin{align}
Q_{16}(b) =  &b^{16}+128.276 b^{15}+7359.58 b^{14}+249796. b^{13}+5.59101\times 10^6 b^{12} \nonumber \\
&+8.70867\times 10^7 b^{11}+9.7108\times 10^8 b^{10}+7.8558\times 10^9 b^9+4.62396\times 10^{10} b^8 \nonumber\\
&+1.9692\times 10^{11} b^7+5.98673\times 10^{11} b^6+1.27226\times 10^{12} b^5
+1.83727\times 10^{12} b^4 \nonumber \\
&+1.74238\times 10^{12} b^3+1.04129\times 10^{12} b^2+3.62215\times 10^{11} b+5.49926\times 10^{10}
\end{align}
The zeros of the polynomial are given by:
\begin{align}
&b_1 =-21.6904 \nonumber\\
&b_2=-18.6188 \nonumber\\
&b_3=-16.1563 \nonumber \\
&b_4=-14.0238\nonumber \\
&b_5=-12.115\nonumber \\
&b_6=-10.3756\nonumber \\
&b_7=-8.77408\nonumber \\
&b_8=-7.29079\nonumber \\
&b_9=-5.91364\nonumber \\
&b_{10}=-4.63575\nonumber \\
&b_{11}=-3.45447\nonumber \\
&b_{12}=-2.37033\nonumber \\
&b_{13}=-1.37913\nonumber \\
&b_{14}=-0.506603-0.513116 i \nonumber \\
&b_{15}=-0.506603+0.513116 i \nonumber \\
&b_{16}=-0.464347
\end{align}
Fitting the first two zeros to the Riemann zeta zeros of the form  $z_i=Ab_i+c$ we find $A=1.52532$ and $c=42.3072$ then the first three zeros determined by the polynomial zeros are $9.22238, 13.90755, 17.6636$ 
which can be compared with the exact first three zeros of the Ramanujan L function which are given by $9.22238, 13.90755, 17.442777$.
One can obtain better estimates for the zeros by going to higher values of $N$ and $p$. It would be interesting to study the behavior of the complex zeros as for large $N$ to see if they go to infinity under the double scaling limit. A plot of the Ramanujan Xi function which is the Baker Akheizer function  for a two matrix model is shown in figure 3.


\begin{figure}
\centering
  \includegraphics[width = .75 \linewidth]{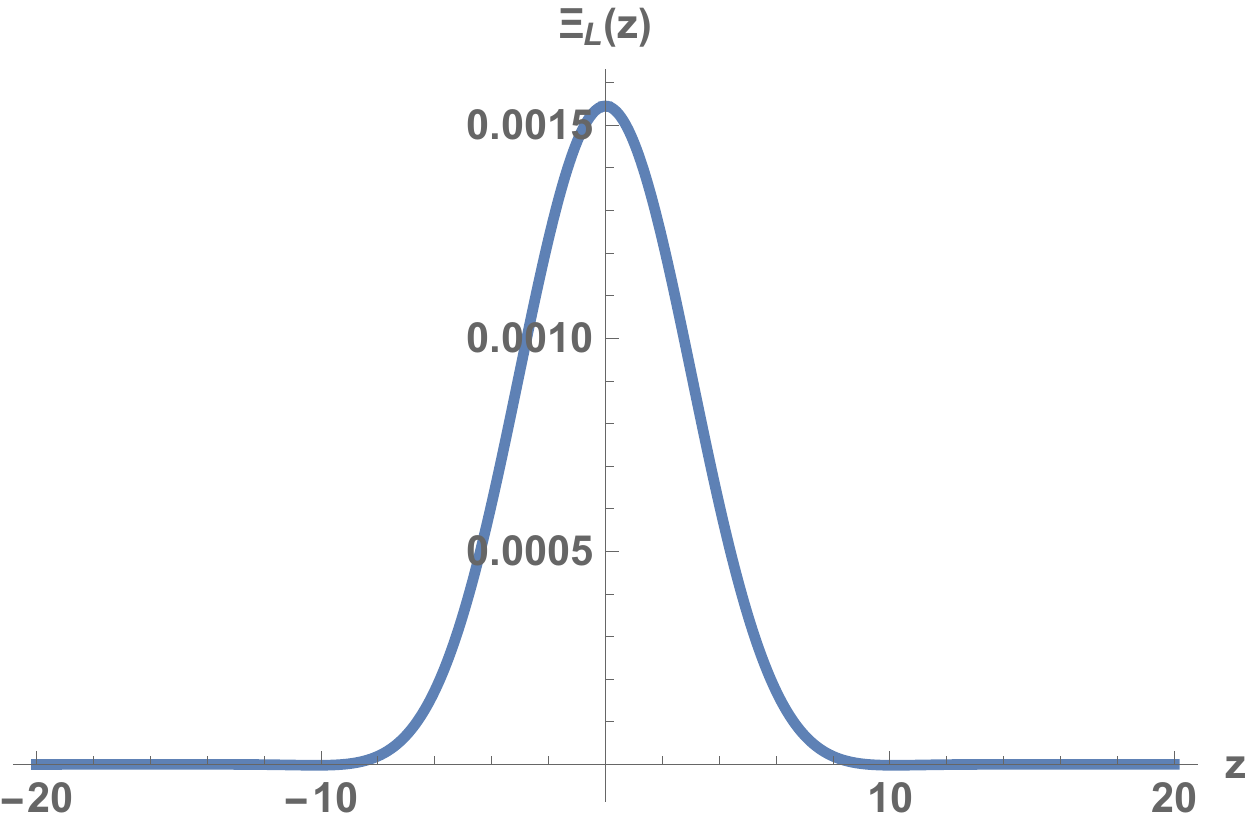}
  \caption{Plot of the Ramanujan Xi function which is the Baker Akheizer function  for a two matrix model associated with the Ramnaujan L function. }
  \label{fig:Radion Potential}
\end{figure}




\section{Two Matrix Model and  Baker-Akhiezer  functions with zeros on the critical line }

In \cite{Paris}
\cite{Bruijn}
\cite{Senouf}
\cite{CardonF}
\cite{Kamimoto}
\cite{Polya}
\cite{Duits6}
\cite{Dominici}
\cite{Ki}
\cite{Dimitrov}
\cite{2013arXiv1309.7019F}
\cite{Cardon2}
\cite{Omar}
\cite{Newman} several criteria were given that lead to Baker-Akhiezer functions with real zeros. If we express the Baker-Akhiezer function as a Fourier integral as:
\begin{equation}\psi (z) = \int_{ - \infty }^\infty  {{e^{ - U(x)}}{e^{izx}}dx} \end{equation}
then the the Baker-Akhiezer function will have real zeros if $U(x)$ is of the form:
\begin{equation}U(x)= \frac{x^{2n}}{2n}\end{equation}
for $n > 2$. 
Another case where Baker-Akhiezer functions with real zeros is when $U(x)$ is of the form: 
\begin{equation}U(x) = \frac{{{x^{4m}}}}{{2m}} + s_{2m-1}\frac{{{x^{2m}}}}{{2m}} + s_1\frac{{{x^2}}}{2}\end{equation}
with $s_{2m-1},s_1$ real $m = 1,2,3, \dots $,2 with $s_1 \le 0$.  
Finally if $U(x)$ has a leading term $\frac{x^{2n}}{2n}$ and it's value at imaginary arguments $U(iu)$ has real zeros only than the Baker-Akhiezer function will have real zeros.

An example of the first type above where $U(p)= \frac{x^{2n}}{2n}$  with $n=4$ and the Baker-Akhiezer function given by the generalized Airy function. Using formula (3.5) with $N=16$ we have the polynomial:
\begin{align}
Q_{16}(b)= &b^{16}+112. b^{15}+5722.5 b^{14}+176798. b^{13}+3.69077\times 10^6 b^{12} \nonumber \\
&+5.51006\times 10^7 b^{11}+6.07273\times 10^8 b^{10}
+5.02779\times 10^9 b^9+3.15136\times 10^{10} b^8 \nonumber \\
&+1.49533\times 10^{11} b^7+5.33211\times 10^{11} b^6+1.40685\times 10^{12} b^5+2.67577\times 10^{12} b^4 \nonumber \\
&+3.51786\times 10^{12} b^3+2.98205\times 10^{12} b^2+1.43339\times 10^{12} b+2.87152\times 10^{11}
\end{align}
with zeros give by:
\begin{align}    
&b_1 = -15.5108\nonumber \\
&b_2 = -13.7506\nonumber \\ 
&b_3 = -12.3194\nonumber \\ 
&b_4 = -11.0623\nonumber \\
&b_5 = -9.91975\nonumber \\ 
&b_6 = -8.86087\nonumber \\ 
&b_7 = -7.86714\nonumber \\ 
&b_8 = -6.92631\nonumber \\ 
&b_9 = -6.02971\nonumber \\ 
&b_{10} = -5.17091\nonumber \\ 
&b_{11} = -4.34496\nonumber \\ 
&b_{12} = -3.54788\nonumber \\
&b_{13} = -2.77634\nonumber \\ 
&b_{14} = -2.02742\nonumber \\ 
&b_{15} = -1.29834\nonumber \\ 
&b_{16} = -0.587221\nonumber \\
\end{align}
Transforming these eigenvalues of the form $y=Ab+c$ these eigenvalues with $A = 1.43304,c= 24.7926$ we have the first three zeros as:
\begin{align}
y_1 = 2.56503 \nonumber \\
y_2 = 5.08746 \nonumber \\
y_3 = 7.13834
\end{align}
defining:
\begin{equation}Ai_{(7,1)}\left( {y,{s_1},{s_3},{s_5}} \right) = \int_{ - \infty }^\infty  {dx{e^{ - \left( {{s_1}\frac{{{x^2}}}{2} + {s_3}\frac{{{x^4}}}{4} + {s_5}\frac{{{x^6}}}{6} + \frac{{{x^8}}}{8}} \right)}}} {e^{ixy}}\end{equation}
which can be compared to the first three zeros of $Ai_{(7,1)}(y)=Ai_{(7,1)}(y,0,0,0)$ given by $2.56503, 5.08746, 7.53357$. A plot of $Ai_{(7,1)}(y)$  showing the first few zeros is in in figure 4.
\begin{figure}
\centering
  \includegraphics[width = .75 \linewidth]{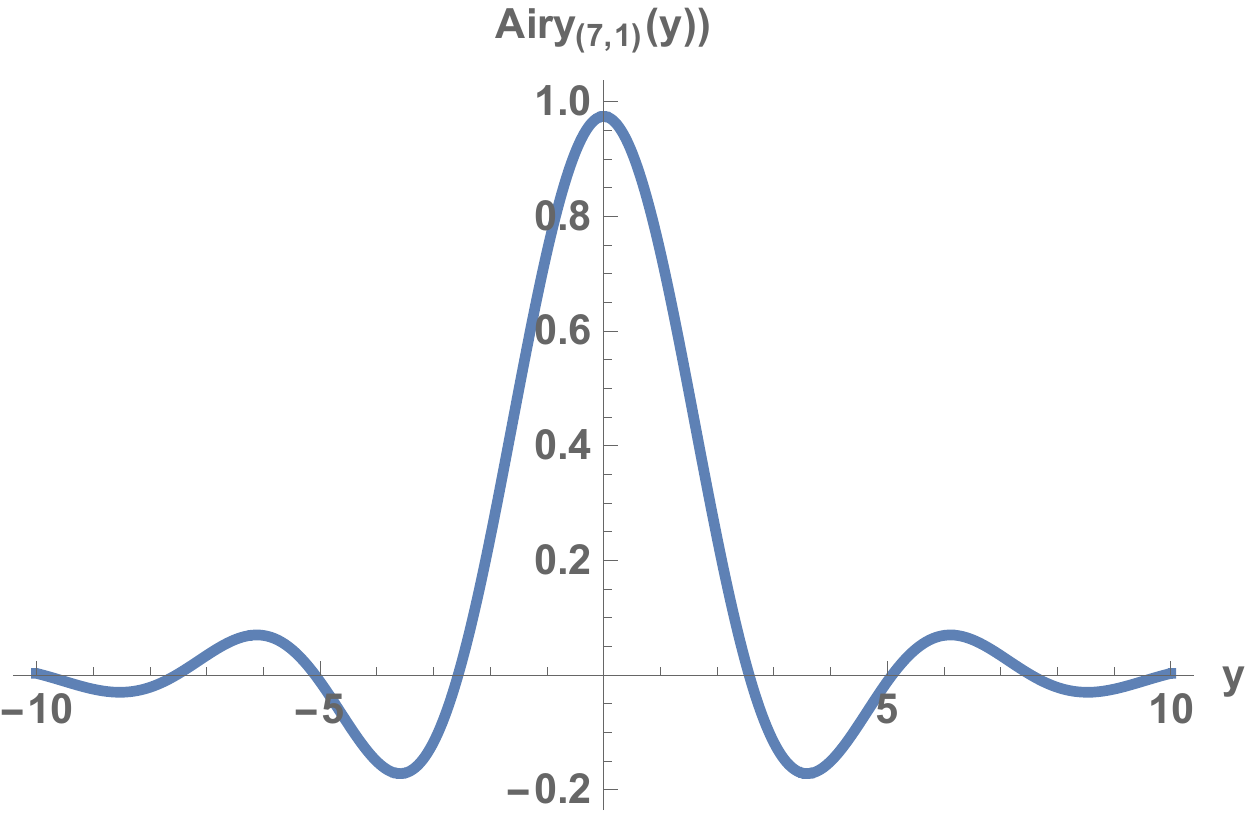}
  \caption{Plot of the generalized Airy function which is the Baker Akheizer function  for the $(7,1)$ two matrix model showing the first few zeros. }
  \label{fig:Radion Potential}
\end{figure}

For an example of the second  type of $U(x)$ that yield a Baker-Akhiezer function that will have real zeros we will take:
\begin{equation}U(x) = \frac{{{x^8}}}{8} + 3\frac{{{x^4}}}{4} - \frac{{{x^2}}}{2}\end{equation}
so for this case:
\begin{align}
    &s_1=-1\nonumber \\
    &s_3 = 3 \nonumber \\
    &s_5= 0 
\end{align}
and for $N=16$ the polynomial representation of the expectation value of the characteristic polynomial is:
\begin{align}
Q_{16}(b)=&b^{16}+146. b^{15}+9708.52 b^{14}+389701. b^{13}+1.05503\times 10^7 b^{12} \nonumber\\
&+2.03871\times 10^8 b^{11}+2.9022\times 10^9 b^{10}+3.09666\times 10^{10} b^9+2.49541\times 10^{11} b^8 \nonumber \\
&+1.51836\times 10^{12} b^7+6.92309\times 10^{12} b^6+2.32847\times 10^{13} b^5+5.62629\times 10^{13} b^4 \nonumber \\
&+9.36259\times 10^{13} b^3+1.00052\times 10^{14} b^2+6.03675\times 10^{13} b+1.51162\times 10^{13}
\end{align}
whose zeros are:
\begin{align}
&b = -20.6811 \nonumber \\
&b = -18.2512 \nonumber \\
&b = -16.2823  \nonumber \\
&b = -14.5587  \nonumber \\
&b = -12.9977 \nonumber \\
&b = -11.5564 \nonumber \\
&b = -10.2094 \nonumber \\
&b = -8.93992  \nonumber \\
&b = -7.73641 \nonumber \\
&b = -6.59062  \nonumber \\
&b = -5.49659  \nonumber \\
&b = -4.45016 \nonumber \\
&b = -3.44881  \nonumber \\
&b = -2.49194 \nonumber \\
&b = -1.58201 \nonumber \\
&b = -0.726718 
\end{align}
Fitting the first two zeros to the generalized Airy zeros of the form  $z_i=Ab_i+c$ we find $A=1.27473$ and $c=29.2615$ then the first three zeros determines by the polynomial zeros are $ 2.89881, 5.99627 , 8.50607$ which can be compared with the exact values $2.89881, 5.99627 , 8.6996 $ for the first three zeros of $Ai_{(7,1)}(y,-1,3,0)$. One can increase the accuracy of the polynomial zero estimation of the zeros by going to higher values of $N$.

For an example of the third  type of $U(x)$ that yield a Baker-Akhiezer function will have real zeros we will take:
\begin{equation}U(x) = \frac{{{x^8}}}{8} + 3\frac{{{x^6}}}{6} + 3\frac{{{x^4}}}{4} + \frac{{{x^2}}}{2}\end{equation}
so that $U'(iu) = 0$ has solutions that are real.
For this case:
\begin{align}
    &s_1=1\nonumber \\
    &s_3 = 3 \nonumber \\
    &s_5= 3 
\end{align}
and for $N=16$ the polynomial representation of the expectation value of the characteristic polynomial is:

\begin{align}
Q_{16}(b)=&b^{16}+270. b^{15}+33203.7 b^{14}+2.46483\times 10^6 b^{13}+1.23403\times 10^8 b^{12}  \nonumber \\
&+4.40938\times 10^9 b^{11}+1.1605\times 10^{11} b^{10}+2.28877\times 10^{12} b^9+3.40782\times 10^{13} b^8 \nonumber \\
&+3.82914\times 10^{14} b^7+3.22164\times 10^{15} b^6+1.99709\times 10^{16} b^5+8.87895\times 10^{16} b^4 \nonumber \\
&+2.71163\times 10^{17} b^3+5.29656\times 10^{17} b^2+5.80173\times 10^{17} b+2.60471\times 10^{17}
\end{align}
whose zeros are given by:
\begin{align}
&b_1 = -38.1637 \nonumber \\
&b_2 = -33.7063 \nonumber \\
&b_3 = -30.0921 \nonumber \\
&b_4 = -26.9256 \nonumber \\
&b_5 = -24.0553 \nonumber \\
&b_6 = -21.4028  \nonumber \\
&b_7 = -18.9211 \nonumber \\
&b_8 = -16.5794 \nonumber \\
&b_9 = -14.3561 \nonumber \\
&b_{10} = -12.2356  \nonumber \\
&b_{11} = -10.2064 \nonumber \\
&b_{12} = -8.26003  \nonumber \\
&b_{13} = -6.39036 \nonumber \\
&b_{14} = -4.59403  \nonumber \\
&b_{15} = -2.87196 \nonumber \\
&b_{16} = -1.23928 
\end{align}
Fitting the first two zeros to the generalized Airy zeros of the form  $z_i=Ab_i+c$ we find $A=0.790259$ and $c=34.3341$ then the first three zeros determines by the polynomial zeros are $ 4.17486, 7.69736, 10.5535$ which can be compared with the exact values $ 4.17486, 7.69736, 10.9217 $ for the first three zeros of $Ai_{(7,1)}(y,1,3,3)$. One can increase the accuracy of the polynomial zero estimation of the zeros by going to higher values of $N$.

\section{Two Matrix models and the infinite $p$ limit}

There has been recent interest in the infinite $p$ limit of two Matrix models. In \cite{Mahajan:2021nsd} the $(p,2)$ two matrix model was used in the limit of infinite $p$ to study the $c \rightarrow \infty$ Liouville theory in the semiclassical limit. The coupled $(p,q)$ minimal matter model has central charge:
\begin{equation}c = 1 - 6\frac{{{{\left( {p - q} \right)}^2}}}{{pq}}\end{equation}
which for the $(p,1)$ model yields:
\begin{equation}c = 1 - 6\left( {p - 2 + \frac{1}{p}} \right)\end{equation}
which goes to $-\infty$ as $p$ goes to infinity. So the Liouvlle threory has central charge 
$+\infty$ so that the total theory plus ghosts is conformally invariant. The large $p$ limit allows one to use the semiclassical limit for the Liouville theory as in \cite{Mahajan:2021nsd}. The large $p$ limit will also allow us to construct a two Matrix model which more accuractly computes the zeros of the Baker-Akhiezer function. 

For the potential:
\begin{equation}U(x)=\cosh\ (x)\end{equation}
we have the expansion:
\begin{equation}
\cosh (x) = \sum\limits_{n = 0}^\infty  {\frac{{{x^{2n}}}}{{\left( {2n} \right)!}}} \end{equation}
So we can obtain exact expressions for the $s_n$ for the corresponding two matrix model. Only the odd $n$ are nonzero which take the value:
\begin{equation}
{s_{2j + 1}} = \frac{1}{{\left( {2j + 1} \right)!}}{\left( {p!} \right)^{\frac{{2j + 2}}{{p + 1}}}}
\end{equation}
for $2j + 1 = 1,3, \ldots ,p - 2 $.

The Baker-Akhiezer function associated with this $U(x)$ is the K-Bessel function 
\begin{equation}
{\psi(z)=K_{iz}}(1) = \int_{ - \infty }^\infty  {{e^{ - \cosh x}}{e^{izx}}dx} 
\end{equation}
The potential for the two matrix model in this case
\begin{equation}
V(A,B)=Tr(\sum _{j=0}^{\frac{p-3}{2}} \frac{(p!)^{\frac{(2 j+1)+1}{p+1}} \epsilon ^{p-(2 j+1)} V(2 j+1,A+1)}{(2 j+1)!}+V(p,A+1)-AB)
\end{equation}
Because we have a closed form expressing for the potential as for arbitrary $p$ we can consider the large $p$ limit to closer approximate the zeros of  Baker-Akhiezer function using the two matrix model. For $p=7$ and $N=16$ we have:
\begin{align}
&s_1=8.42573 \nonumber \\
&s_3=11.8322 \nonumber \\
&s_5=4.98473
\end{align}
and the $Q_{N}(b)$ polynomial is:
\begin{align}
Q_{16}(b)=&1. b^{16}+470.227 b^{15}+100563. b^{14}+1.29618\times 10^7 b^{13}+1.12484\times 10^9 b^{12} \nonumber \\
&+6.95391\times 10^{10} b^{11}+3.16015\times 10^{12} b^{10}+1.07377\times 10^{14} b^9+2.74766\times 10^{15} b^8 \nonumber \\
&+5.29135\times 10^{16} b^7+7.60608\times 10^{17} b^6+8.02664\times 10^{18} b^5+6.04942\times 10^{19} b^4 \nonumber \\
&+3.11591\times 10^{20} b^3+1.01997\times 10^{21} b^2+1.85683\times 10^{21} b+1.36908\times 10^{21}
\end{align}
with zeros given by:
\begin{align}
&b_1 = -67.7171 \nonumber \\
&b_2= -59.5952 \nonumber \\
&b_3= -53.0264 \nonumber \\
&b_4 = -47.2861 \nonumber \\
&b_5 = -42.0962 \nonumber \\
&b_6 = -37.3136 \nonumber \\
&b_7 = -32.8525 \nonumber \\
&b_8 = -28.657 \nonumber \\
&b_9 = -24.6889 \nonumber \\
&b_{10} = -20.9208 \nonumber \\
&b_{11} = -17.3337 \nonumber \\
&b_{12} = -13.9147 \nonumber \\
&b_{13} = -10.6571 \nonumber \\
&b_{14} = -7.56181 \nonumber \\
&b_{15} = -4.64341 \nonumber \\
&b_{16} = -1.96262
\end{align}
We can relate then zeros to the zeros of $K_{iz}(1)$ through a linear relation $z_i=Ab_i +c$ with $A= 0.193542$ and $c=16.0687$. This can be used to obtain estimates the first three zeros as $2.96255, 4.53449, 5.80583$ with the first two fixed by the linear relation. This can be compared with the exact first three zeros of $K_{iz}(1)$ which are $2.96255, 4.53449, 5.87987$. A plot of  $K_{iz}(1)$ showing the first few zeros is in figure 5.
\begin{figure}
\centering
  \includegraphics[width = .75 \linewidth]{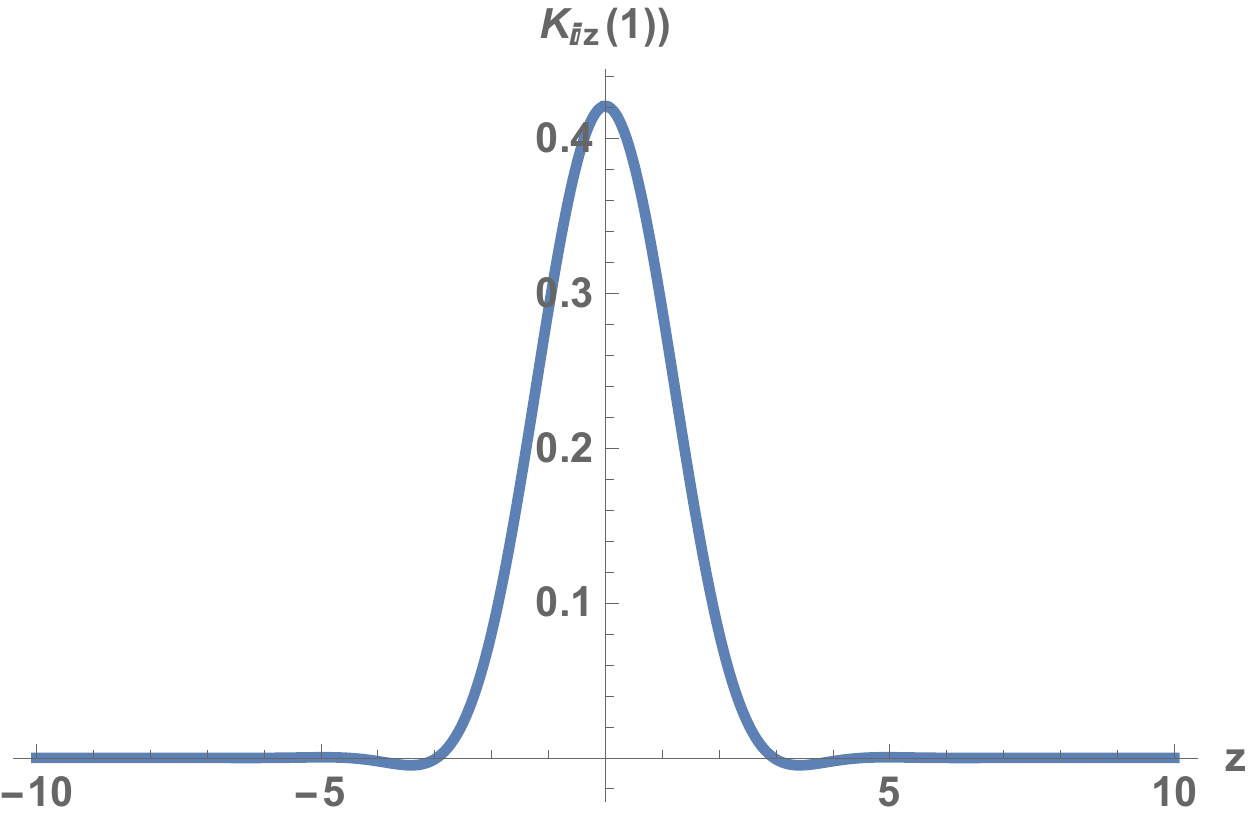}
  \caption{Plot of the Bessel $K_{iz}(1)$  function which is the Baker Akheizer function  for the $(\infty,1)$ two matrix model. }
  \label{fig:Radion Potential}
\end{figure}
It is also interesting in this case that
\begin{equation}
U'(iu) = \sin (u)
\end{equation}
only has real zeros so the third criteria for the Baker-Akhiezer function to have real zeros holds for this case too. Notably this seems to be true for the $U(x)$ functions associated with the Riemann Xi function and Ramanujan L-function. Those cases additionally have modular invariance properties. The function $U'(x)$  has an additional interpretatation as a superpotential whose modular properties are similar to those found in supersymmetric string compactifications.

\subsection{Large $p$ and zeta zeros}

To study the the zeta zeros at large $p$ it is better to use the product of the gamma function and Dirichlet eta function which has a more strightforward expansion for the function $U(x)$.  The product of the gamma function and Dirichlet eta function share nontrivial zeta zeros with the Riemann zeta function as shown in figure 6. The function $U(x)$  is defined through the integral representation as :
\[\Gamma \left( {iz + \frac{1}{2}} \right)\eta \left( {iz + \frac{1}{2}} \right) = \int_{ - \infty }^\infty  {{e^{ - U(x)}}{e^{izx}}dx} \]
with:
\begin{equation}
U(x)=-\log \left(\frac{e^{-\frac{1}{2} (x+\log (2))}}{\exp \left(e^{-(x+\log (2))}+1\right)}\right)
\end{equation}
Expanding $U(x)$  out to $x^{20}$ and rescaling $x$ so the the $x^{20}$ term has coefficient $1/20$ we have:
\begin{align}
U(x)= &0.05 x^{20}-0.135112 x^{19}+0.34685 x^{18}-0.843543 x^{17}+1.93754 x^{16}-4.18855 x^{15}+8.48885 x^{14} \nonumber \\
&-16.0572 x^{13}+28.2038 x^{12}-45.7281 x^{11}+67.9625 x^{10}-91.8254 x^9+111.66 x^8-120.693 x^7 \nonumber \\
&+114.15 x^6-92.538 x^5+62.5149 x^4-33.7861 x^3+13.6947 x^2
\end{align}
so we have $s_n$ deformation parameters for the corresponding $(19,1)$ two matrix model:
\begin{align}
&s_1=13.6947\nonumber \\
&s_2=-33.7861\nonumber \\
&s_3=62.5149\nonumber \\
&s_4=-92.538\nonumber \\ 
&s_5=114.15\nonumber \\
&s_6=-120.693\nonumber \\
&s_7=111.66\nonumber \\
&s_8=-91.8254\nonumber \\
&s_9=67.9625\nonumber \\
&s_{10}=-45.7281\nonumber \\
&s_{11}=28.2038\nonumber \\
&s_{12}=-16.0572\nonumber \\
&s_{13}=8.48885\nonumber \\
&s_{14}=-4.18855\nonumber \\
&s_{15}=1.93754\nonumber \\
&s_{16}=-0.843543\nonumber \\
&s_{17}=0.34685\nonumber \\
&s_{18}=-0.135112 
\end{align}
Then using formula (3.5) and (3.9) we have for $Q_{16}(b)$:
\begin{align}
Q_{16}(b)=& b^{16}+95.0376 b^{15}+3997.72 b^{14}+98185.2 b^{13} \nonumber \\
&+1.56353\times 10^6 b^{12}
+1.69426\times 10^7 b^{11}+1.27431\times 10^8 b^{10}+6.64962\times 10^8 b^9 \nonumber \\
&+2.35587\times 10^9 b^8+5.36099\times 10^9
b^7+6.78531\times 10^9 b^6+2.27449\times 10^9 b^5 \nonumber \\
&-4.40915\times 10^9 b^4-3.60649\times 10^9 b^3+1.09591\times 10^9 b^2+1.06809\times 10^9 b-3.43796\times 10^7
\end{align}
with zeros given by:
\begin{align}
&b_1 = -17.0397\nonumber \\ 
&b_2 = -14.5462\nonumber \\ 
&b_3 = -12.5532\nonumber \\
&b_4 = -10.8316\nonumber \\ 
&b_5 = -9.29383\nonumber \\
&b_6 = -7.89499\nonumber \\
&b_7 = -6.60873\nonumber \\
&b_8 = -5.41845\nonumber \\ 
&b_9 = -4.3135\nonumber \\
&b_{10} = -3.28741\nonumber \\
&b_{11} = -2.33716\nonumber \\
&b_{12} = -1.46324\nonumber \\
&b_{13} = -0.670506\nonumber \\
&b _{14}= 0.0312907\nonumber \\
&b_{15} = 0.594787 - 0.166798 i\nonumber \\
&b_{16} = 0.594787 + 0.166798 i\nonumber \\
\end{align}
We can relate then zeros to the zeros of $\Gamma \left( {iz + \frac{1}{2}} \right)\eta \left( {iz + \frac{1}{2}} \right)$ through a linear relation $z_i=Ab_i +c$ with $A= 2.7621$ and $c=61.2001$. This can be used to obtain estimates the first three zeros as $14.1347, 21.022, 26.527 $ with the first two fixed by the linear relation. This can be compared with the exact first three zeros of $\Gamma \left( {iz + \frac{1}{2}} \right)\eta \left( {iz + \frac{1}{2}} \right)$ which are $ 14.1347, 21.022, 25.0109$. So we see that even going to large $p$ using the $(19,1)$ two matrix model we still have zeros off the critical line,   however again these zeros may go to infinity as $N$ goes to infinity.
\begin{figure}
\centering
  \includegraphics[width = .75 \linewidth]{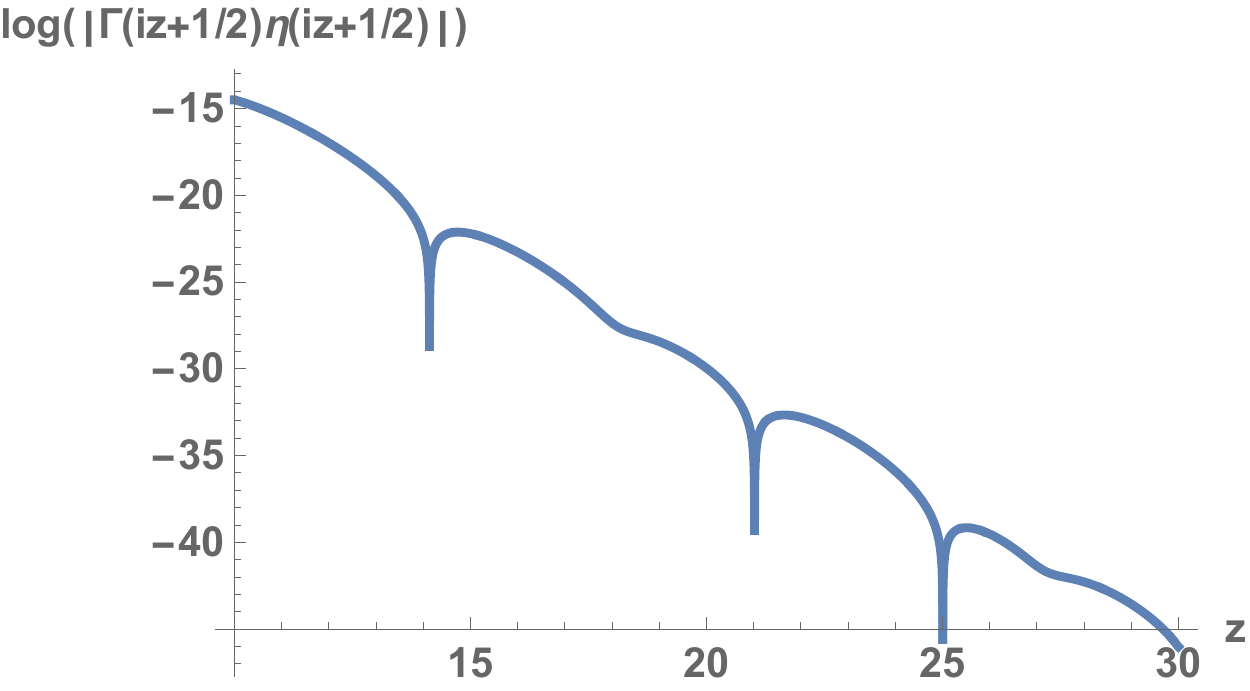}
  \caption{Plot of the logarithm of the absolute vale of $\Gamma \left( {iz + \frac{1}{2}} \right)\eta \left( {iz + \frac{1}{2}} \right)$ function which is the Baker Akheizer function  for the $(\infty,1)$ two matrix model. }
  \label{fig:Radion Potential}
\end{figure}
\section{Obstructions to the Master Matrix for the Two Matrix Model}

The fact that the Baker-Akhiezer function for the two matrix model can have complex zeros can lead to an obstruction to the construction of the master matrix for the two matrix model. This is because the master matrix is a special Hermitian matrix whose characteristic polynomial yields the expectation value for the two matrix model in the limit of infinite $N$. That is:
\begin{equation}
{\left\langle {\det (b - B} \right\rangle _{N \times N}} = \det (b - {B_{master}})
\end{equation}
Since $B_{master}$ is Hermitian it's characteristic polynomial must have real zeros. But we have seen the exact expression for the expectation of the characteristic polynomial in terms of the biorthogonal polynomial $Q_N(b)$ 
which can have complex zeros depending on the potential of the two matrix model. In the large $N$ limit these complex zeros can lead to complex zeros of the Baker-Akhiezer function. The master matrix can never lead to these complex zeros thus we say that there is an obstruction to the construction of the master matrix in those cases.

 There are several approaches to construct the master field or master matrix \cite{Witten:1979pi}
\cite{Coleman:1985rnk}
\cite{tHooft:1973alw}
\cite{Carlson:1982vz}
\cite{Halpern:1984ty}
\cite{Halpern:1999ra}
\cite{Haan:1981ks}
\cite{Gopakumar:2002wx}
\cite{Haan:1980cw}
\cite{Accardi:1995ju}
\cite{Engelhardt:1996da}
\cite{Kuroki:1998rx}
\cite{Douglas:1994kw}. One has the non-commutative probability algebraic approach of \cite{Gopakumar:1994iq}. Here the master matrix take the form of a combination the raising and lowering operators. For example for the $(2,1)$ matrix model
the master matrix or operator is $a + {a^\dag }$ with creation and annihilation operators which leads to the representation of the characteristic polynomial in terms of the Hermite polynomials through the relation to the Harmonic oscillator. As discussed above these leads to the Airy function as the Baker-Akhiezer function in the double scaling limit. Notably the zeros of the Baker-Akhiezer function are real which is consistent with the existence of a master matrix in that case. In other cases the polynomials for the two matrix model are not related  to a standard quantum mechanical system like the harmonic oscillator. This is especially  true for two matrix models whose orthogonal polynomials have complex zeros.  

Another  approach to the master matrix is the quenched master field \cite{Greensite:1982mf}
\cite{Greensite:1982ck}\cite{Alberty:1983ti}\cite{Klinkhamer:2021wrv}
which can be applied to the matrix models discussed in the previous sections. Here one defines the partition function:
\begin{equation}
Z(g) = \int {dAdB{e^{ - S(A,B)}}} 
\end{equation}
where $S(A,B)$  is the action of the two matrix model given by:
\begin{equation}
S(A,B) = \frac{1}{g}Tr(V(A + 1) - AB)
\end{equation}
One then introduces a stochastic time $\tau$ and a diagonal matrix $\hat p$ as well as stochastic time independent matrices $\hat a$ and $\hat b$ through:
\[\hat A = {e^{i \hat p\tau }}\hat a{e^{ - i \hat p\tau }}\]
\begin{equation}\hat B = {e^{i \hat p\tau }}\hat b{e^{ - i \hat p\tau }}\end{equation}
The quenched master field equations are then given by:
\[i\left( {{{\hat p}_k} - {{\hat p}_\ell }} \right){{\hat a}_{k\ell }} =  - {\left. {\frac{{\delta S}}{{\delta {A_{\ell k}}}}} \right|_{A = \hat a,B = \hat b}} + \hat \eta _{k\ell }^{(1)}\]
\begin{equation}i\left( {{{\hat p}_k} - {{\hat p}_\ell }} \right){{\hat b}_{k\ell }} =  - {\left. {\frac{{\delta S}}{{\delta {B_{\ell k}}}}} \right|_{A = \hat a,B = \hat b}} + \hat \eta _{k\ell }^{(2)}\end{equation}
which reduce to:
\begin{align}
&i\left( {{{\hat p}_k} - {{\hat p}_\ell }} \right){{\hat a}_{k\ell }} =  - \frac{1}{g}V'(a + I)_{k \ell} + \frac{1}{g}{{\hat b}_{k\ell }} + \hat \eta _{k\ell }^{(1)}\nonumber \\
&i\left( {{{\hat p}_k} - {{\hat p}_\ell }} \right){{\hat b}_{k\ell }} =   \frac{1}{g}{{\hat a}_{k\ell }} + \hat \eta _{k\ell }^{(2)}\end{align}
In these equations $\hat p_i$ are  randomly distributed master momenta and  $\hat \eta _{k\ell }^{(2)}$  are master noise matrices whose components are Gaussian distributed random numbers.
 The solution to the quenched master field  equations can be formulated as an optimization problem by forming:
\begin{align}
    & E_{k\ell}=i\left( {{{\hat p}_k} - {{\hat p}_\ell }} \right){{\hat a}_{k\ell }}   + \frac{1}{g}V'(a + I)_{k\ell} - \frac{1}{g}{{\hat b}_{k\ell }} - \hat \eta _{k\ell }^{(1)}\nonumber \\
& F_{k \ell} = i\left( {{{\hat p}_k} - {{\hat p}_\ell }} \right){{\hat b}_{k\ell }} -   \frac{1}{g}{{\hat a}_{k\ell }} - \hat \eta _{k\ell }^{(2)}\end{align}
and optimizing:
\begin{equation}
C = \sum\limits_{i,\ell } {\left( {E_{i\ell }^2 + F_{i\ell }^2} \right)} 
\end{equation}
Various optimization procedures can be used for finding the minimum of $C$. It is interesting that quantum computing can be used for optimization problems so this is another area of potential application of quantum algorithms to the Riemann hypothesis \cite{McGuigan:2023rrn}
\cite{vanDam}. If the minimum of $C$ is zero than all the equations are satisfied and the value of $ \hat a_{k\ell }$ and $ \hat b_{k\ell }$ determine the Master matrix with real zeros for the expectation value for the characteristic polynomial as $b_{master}$ is Hermitian. If the minimum is not zero than some of the equations could not be satisfied and the master matrix could not be constructed and we say that there is an obstruction. These obstruction could happen if the characteristic polynomial has complex zeros or also if it has real zeros so that an obstruction to the master matrix does not falsify the Riemann hypothesis,  but on the other hand if the master matrix can be constructed the the hypothesis would be true. 

An alternative to the  master matrix equations are the saddle point equations for the two matrix model \cite{Kazakov:2000aq}.
The saddle point approach looks at the eigenvalue integral representation of the partition function:
\begin{equation}Z = \int {\prod\limits_{i = 1}^N {d{a_i}d{b_i}\Delta (a)\Delta (b){e^{\frac{1}{g}\sum\limits_{i = 1}^N {\left( { - V({a_i}} \right) - {b_i} + {a_i}{b_i})} }}} } \end{equation}
with Vandermonde determinants:
\[\Delta (a) = \prod\limits_{i > j} {\left( {{a_i} - {a_j}} \right)} \]
\begin{equation}\Delta (b) = \prod\limits_{i > j} {\left( {{b_i} - {b_j}} \right)} \end{equation}
One then takes the saddle point approximation to the effective potential of the eigenvalues yielding the equations.
\begin{align}
- \frac{1}{g}V'(1 + {a_i}) + \frac{1}{g}{b_i} + &\sum\limits_{j \neq i} {\frac{1}{{{a_i} - {a_j}}}}  = 0\nonumber \\
\frac{1}{g}{a_i} + &\sum\limits_{j \neq i } {\frac{1}{{{b_i} - {b_j}}}}  = 0
\end{align}
It would be interesting to apply it to the Riemann hypothesis in the two matrix model context as the picture of repelling eigenvalues with supporting numerical data about the Riemann zeros. A description in terms of many fermion ground state wave functions is also possible with the Vandermonde factor being related to a Slater determinant \cite{Callaway:1990nm}
\cite{Cappelli:1993ea}
\cite{Cappelli:2004xk}. Finally new approaches based on positivity and the conformal bootstrap would be interesting to apply to the two matrix model approach to the Riemann hypothesis \cite{Kazakov:2021lel}
\cite{Lin:2023owt}
\cite{Lin:2020mme}
\cite{Han:2020bkb}.

\section{Conclusion}

\begin{table}[h]
\centering
\begin{tabular}{|l|l|l|l|l|l|l|}
\hline
Baker-Akhiezer function      & $U(x)$  & $z_3  (N=16)$ &  $z_3$ exact &  On CL finite N & infinite N \\ \hline
$Ai(z)$   &  $ i \frac{x^3}{3} $ &  -5.56709 & -5.52056 & Y  &  Y  \\ \hline
Riemann  $\Xi(z)$    & $ -\log \left(\Phi (x)\right) $  & 26.5505 & 25.0109 & N & ? \\ \hline
Ramanujan  $\Xi_L(z)$  & $ -\log \left(\Phi_L (x)\right) $  & 17.6636 &  17.442777 & N & ? \\ \hline
$Ai_{(7,1)}(z)$ &  $\frac{x^8}{8} $  & 7.13834 & 7.53357 & Y & Y \\ \hline
$Ai_{(7,1)}(z,-1,3,0)$  & $\frac{x^8}{8} + 3\frac{x^4}{4} - \frac{x^2}{2}$ & 8.50607 & 8.6996 & Y & Y\\ \hline
$Ai_{(7,1)}(z,1,3,3)$ & $\frac{x^8}{8} + 3\frac{x^6}{6}+ 3\frac{x^4}{4} + \frac{x^2}{2}$ & 10.5535 & 10.9217 & Y & Y  \\ \hline
$K_{i z}(1) $  & $ \cosh(x) $  & 5.80583 & 5.87987 & Y & Y \\ \hline
$\Gamma \left( {iz + \frac{1}{2}} \right)\eta \left( {iz + \frac{1}{2}} \right)$  & $ \log \left(\frac{e^{-\frac{1}{2} (x+\log (2))}}{\exp \left(e^{-(x+\log (2))}+1\right)}\right) $  & 26.527 & 25.0109 & N & ? \\ \hline
\end{tabular}
\caption{\label{tab:table-name} Results from studying the Baker-Akhiezer function for $(p,1)$  two matrix models. The third zero of the Baker-Akheizer function was estimated using the polynomial approach for $N=16$ and compared with the exact value. The first two zeros were fixed by a linear transformation. The only Baker-Akhiezer functions we studied with zeros off the critical line (CL) at finite $N$ were the Riemann Xi, Ramanujan L function and the product of gamma and eta functions. Whether these functions have zeros off the critical line as $N$ goes to infinity is unknown and its solution is equivalent to the Riemann and generalized Riemann hypothesis.} 
\end{table}

In this paper we studied the Riemann Xi function and other functions as Baker-Akhiezer functions for two matrix models. We recorded our results in table 1. For finite $N$ we  found zeros off the critical line for the Riemann Xi and Ramanujan L functions using the method of bi-orthogonal polynomials applied to the two matrix model. We studied $N=16$ and it would be interesting to study larger $N$ to see what happens to the zeros off the critical line as $N$ goes to infinity. It would also be interesting to apply other methods to the two matrix model associated with the Riemann Xi function such as the master matrix approach, saddle point method or bootstrap approach. This may give additional information on the fate of the zeros off the critical line at large $N$. Other methods that can be applied include unitary matrix models, phase space and prime factorization applied to the Riemann zeros discussed in \cite{Dutta:2016byx}
\cite{Chattopadhyay:2017ckc}
\cite{Chattopadhyay:2018bzs}.

For both the one and two matrix model the characteristic polynomial of the Jacobi matrix yields the the expectation value of the characteristic polynomial which becomes the Baker-Akhiezer function in the large $N$ limit.  Unlike the one matrix model, for the two matrix model the Jacobi matrix is not symmetric and for a range of coupling $g$  and deformation parameters $s_i$ can develop complex eigenvalues for finite $N$. The important question is what happens to these complex eigenvalues that potentially become zeros off the critical in the scaling limit as $N$ goes to infinity. It is thought that solution to the Riemann hypothesis must in involve the arithmetic properties of the zeta function. The modular invariance properties of the potentials $-\log(\Phi(x))$ and $-\log (\Phi_L(x))$ functions may also play a role as these functions can be expressed in terms of theta functions which also have arithmetic properties. Interestingly modular properties of potentials play a role in fundamental physics  through modular invariant axion-dilaton potentials \cite{Cvetic:1993ty}
\cite{Garcia-Bellido:1992its}
\cite{Horne:1994mi}
\cite{Gonzalo:2018guu}
\cite{Font:1990gx}\cite{Gendler:2022qof}
\cite{Donagi:1996yf}
\cite{Curio:1997rn}\cite{Leedom:2022zdm}
\cite{Alexander:2023qym}
\cite{Cribiori:2023sch}. It would be remarkable (although not unprecedented) if a fundamental problem in mathematics such as the Riemann hypothesis turned out to be related to  fundamental problems in the physics of elementary particles.

\section*{Acknowledgements}
We wish to thank Arghya Chattopadhyay
for interesting suggestions and references relating the Riemann Zeta function to Matrix models and Berry-Keating Hamiltonians. We also thank Nicole Righi foor sending references related to modular potentials for the heterotic string and possible metastable de Sitter vacua.

\end{document}